\input amstex
\input xy
\input epsf
\xyoption{all}
\documentstyle{amsppt}
\document
\magnification=1200
\NoBlackBoxes
\nologo
\hoffset1.5cm
\voffset2cm

\def\cO{\Cal{O}}

\def\d{\delta}

\def\cO{\Cal{O}}

\def\mathcal A{\Cal{A}}
\def\ra{\to}

\pageheight {16cm}

\bigskip

\centerline{\bf  ARITHMETIC DIFFERENTIAL EQUATIONS}
\smallskip
\centerline{\bf OF PAINLEV\'E VI TYPE} 

\medskip
\centerline{\bf Alexandru Buium${}^1$, Yuri I. Manin${}^2$}

\medskip

\centerline{\it ${}^1$University of New Mexico, Albuquerque}
\smallskip
\centerline{\it ${}^2$Max--Planck--Institut f\"ur Mathematik, Bonn, Germany}

\bigskip

ABSTRACT. {\it Using the description of Painlev\'e VI family of
differential equations in terms of  a universal elliptic curve,
going back to R.~Fuchs (cf. [Ma96]), we translate it into
the realm of Arithmetic Differential Equations (cf. [Bu05]),
where the role of derivative ``in the  $p$--adic direction''
is played by a version of Fermat quotient.}

\vskip1cm

\centerline{\bf Introduction and brief summary}

\medskip

This article is dedicated to  the study of differential equations of the
Painlev\'e VI type ``with $p$--adic time'', in which the role of
time derivative is played by a version of $p$--Fermat quotient.
Richness of $p$--adic differential geometry was already demonstrated
in many papers of the first--named author: see in particular the monograph [Bu05].

\smallskip

Applicability of this technique to the Painlev\'e equations is ensured
by the combination of  two approaches:  R.~Fuchs's treatment of the classical case modernised in
[Ma96], and constructions of $p$--adic differential characters in [Bu95].

\smallskip

In sec.~1 below we present a short introduction to the $p$--adic differential geometry. In sec.~2 and 3
we introduce several versions of $p$--adic PVI. Sec.~4 is dedicated to the
problem of transposing into $p$--adic domain of main features of Hamiltonian formalism.
We have not found a definitive answer to this problem, and Painlev\'e equations serve here
as a very stimulating testing ground. Finally, section 5 treats another problem, suggested
by comparison of several versions of $p$--adic PVI, but presenting an independent interest.

\smallskip

{\it Acknowledgement.} This article was conceived during the Spring 2013 Trimester program
``Arithmetic and Geometry'' of the Hausdorff Institute for Mathematics (HIM), Bonn. The authors are grateful to HIM for
stimulating atmosphere and working conditions.

\bigskip

\centerline{\bf  1. Arithmetic differential equations: background}

\medskip

 We start with a brief summary of relevant
material from [Bu05], Chapters 2, 3. Fix a prime $p$; in our applications it will be assumed that $p\ge 5$.

\smallskip
{\bf 1.1. $p$--derivations.} Recall that, in the conventional commutative algebra, given a ring $A$ and an $A$--module $N$,
{\it a derivation of $A$ with values in  $N$ } is any map 
$\partial :\,A\to N$ such that the map $A\to A\times N:\, a\mapsto (a,\partial a)$
is a ring homomorphism, where $A\times N$ is endowed with the structure
of commutative ring with componentwise addition, and multiplication $(a,m)\cdot (b,n):= (ab, an+bm)$.
Notice that $\{ 0\}\times N$ is the ideal with square zero in $A\times N$.
\smallskip

Similarly, in arithmetic geometry  {\it a $p$--derivation of $A$ with values in an $A$--algebra $B$,
$f:A\to B$, }
is a map $\delta_p:\,A\to B$ such that the map $A\to B\times B:\, a\mapsto (f(a),\delta_p( a))$
is a ring homomorphism $A\to W_2(B)$ where $W_2(B)$ is the ring of $p$--typical Witt vectors of length 2.
Again, if $pB=\{0\}$, Witt vectors of the form $(0,b)$ form the ideal of square zero.
\smallskip
Explicitly, this means that $\delta_p(1)=0$, and
$$
\delta_p(x+y)=\delta_p(x)+\delta_p(y) +C_p(x,y),
\eqno(1.1)
$$
$$
\delta_p(xy)=f(x)^p\cdot \delta_p(y)+f(y)^p\cdot \delta_p(x) +p\cdot \delta_p(x)\cdot\delta_p(y),
\eqno(1.2)
$$
where
$$
C_p(X,Y):=\frac{X^p+Y^p-(X+Y)^p}{p}\in \bold{Z}[X,Y].
\eqno(1.3)
$$
In particular, this implies that for any $p$--derivation $\delta_p:\,A\to B$ the respective map $\phi_p :\, A\to B$
defined by $\phi_p (a):=f(a)^p+p\delta_p(a)$ is a ring homomorphism
satisfying $\phi_p (x)\equiv f(x)^p\,\roman{mod}\,p$, that is ``a lift of the Frobenius map applied to $f$''.
\smallskip

Conversely, having such a lift of Frobenius, we can uniquely reconstruct the respective derivation $\delta_p$
{\it under the condition that $B$ has no $p$--torsion.}

\smallskip

We will often work with $p$--derivations $A\to A$ with respect to the identity map $A\to A$
and keep $p$ fixed; then it might be kept off the notation. Such a pair $(A,\delta )$ is
 called a $\delta$--ring. Morphisms of $\delta$--rings are algebra morphisms
 compatible with their $p$--derivations.

\smallskip

Moreover, our rings (and more generally, schemes) will be $R$--algebras where $R=W(k)$ 
(ring of infinite $p$--typical Witt vectors)
is the completion of the maximal
unramified extension of $\bold{Z}_p$, with residue field
$k:=$ an algebraic closure of $\bold{Z}/p\bold{Z}$. By $\phi :\, R\to R$ we denote
the automorphism acting as Frobenius $x\mapsto x^p$ on $k$, and by $\delta$
the respective $p$--derivation: $\delta(x)= (\phi (x)-x^p)/p.$ The $R$--algebra structure on a $\delta$--ring
is always assumed to be compatible with this
$p$-derivation.

\medskip

{\bf 1.2. Prolongation sequences and $p$--jet spaces.} In the classical situation invoked in 1.1,
there exists an universal derivation 
$$
d:\,A\to \Omega^1(A)
\eqno(1.4)
$$ 
with values in the $A$--module of
differentials.

\smallskip

For $p$--derivations, (1.4) might be replaced by the following construction
(however, see the Remark 1.2.1 below).
\smallskip

Let $A$ be an $R$--algebra. {\it A prolongation sequence} for $A$
consists of a family of $p$--adically complete $R$--algebras $A^i, i\ge 0$, 
where $A^0=A\,\widehat{}$\quad  is the $p$--adic comple\-tion of $A$, and of maps 
$\varphi_i, \delta_i:\,A^i\to A^{i+1}$ satisfying the following conditions:

\smallskip

{\it a) $\varphi_i$ are ring homomorphisms, each $\delta_i$ is a $p$--derivation  with
respect to $\varphi_i$, compatible with $\delta$ on $R$.
\smallskip

b) $\delta_i\circ \varphi_{i-1}=\varphi_i\circ \delta_{i-1}$ for all $i\ge 1$.}
\smallskip

Prolongation sequences  form a category with evident morphisms, ring homomorphisms
$f_i:A^i\to B^i$ commuting with $\varphi_i$ and $\delta_i$,
and in its subcategory with fixed $A^0$ {\it there exists an initial element,
defined up to unique isomorphism} (cf. [Bu05], Chapter 3). It can be called the
universal prolongation sequence.

\smallskip

In the geometric language, if $X=\roman{Spec}\,A$, the formal spectrum of the $i$--th ring $A^i$
in the universal prolongation sequence  is denoted $J^i(X)$ and called 
{\it the $i$--th $p$--jet space of $X$.} Conversely, $A^i=\Cal{O}(J^i(X))$, the ring of global
functions. 

\smallskip
The geometric morphisms (of  formal schemes over $\bold{Z}$) corresponding to $\phi_i$ are denoted
$\phi^i:\, J^i(X)\to J^0(X)=:X\,\widehat{}$ (formal $p$--adic completion of $X$).

\smallskip

This construction is compatible with localisation so that
it can be applied to the non--necessarily affine shemes: cf. [Bu05], Chapter 3.

\smallskip

{\bf 1.2.1. Remark.} Classically, (1.4) extends to the universal map of $A$ to the differential
graded algebra $\Omega^*(A)$, and there is a superficial similarity of this map
with the one, say of $A$ to the inductive limit of its universal prolongation sequence
in the $p$--adic arithmetics context.

\smallskip

However, the classical differential acts in $\bold{Z}_2$--graded supercommutative
algebras and is an {\it odd} operator with square zero, whereas $\delta_p$ are even.

\smallskip

The differential geometry of smooth schemes {\it  in characteristic $p>0$} suggests a
perspective worth exploring. Namely, the sheaf of differential forms on such a scheme
is endowed with the so called {\it Cartier operator $C$}, which is dual to the
Frobenius operator $F:\,\partial \mapsto \partial^p$ acting upon vector fields.
This operator $C$ is $F^{-1}$--linear. One could consider studying $p$--adic lifts
of the Cartier operator from the closed fibre of  the relevant scheme 
to its $p$--adic completion, following the lead of [Bu05]. For a recent survey
of   $F^{-1}$--linear maps  cf. [BlSch12] and references therein.

\medskip

{\bf 1.3. Flows.} Let now $X$ be a smooth affine scheme over $R=W(k)$. Each element
of  $\Cal{O}(J^r(X))$ induces a function  $f:X( R ) \mapsto R$. Such functions are called 
 $\delta$--functions of order $r$ on $X$, and we may and will identify them
 with
 respective elements of  $\Cal{O}(J^r(X))$. For $r=0$,  $\Cal{O}(J^0(X))= \Cal{O}(X)\,\widehat{} $,
 the $p$--adic completion of  $\Cal{O}(X)$.

 \medskip

{\bf 1.3.1. Definition.} {\it a) A system of arithmetic differential equations of order $r$ on $X$ is a subset $\Cal{E}$ of  $\Cal{O}(J^r(X))$.

\smallskip
b) A solution of $\Cal{E}$ is an $R$--point $P\in X( R )$ such that $f(P)=0$ for all $f\in \Cal{E}$.
The set of solutions of $\Cal{E}$ is denoted $Sol(\Cal{E})\subset X( R )$.

\smallskip

c) A prime integral of $\Cal{E}$ is a function $\Cal{H}\in \Cal{O}(X)\,\widehat{}$ such that $\d(\Cal{H}( P ))=0$ for all 
$P\in Sol(\Cal{E})$.}

\smallskip

We will also denote by $Z^r(\Cal{E})$ the closed formal subscheme of $J^r(X)$ generated by $\Cal{E}$.

\smallskip

Now, let $\delta_X$ be a $p$--derivation of  $\Cal{O}(X)\,\widehat{}$. From the universality of the jet
sequence explained in 1.2, it follows that such derivations are in a bijection with the sections
of the canonical morphism  $J^1(X)\to J^0(X)$. 
\medskip

{\bf 1.3.2. Definition.} {\it  The $\d$--{\it flow} associated to $\d_X$ is the system of arithmetic differential equations
of order 1 which is the ideal in $\Cal{O}(J^1(X))$ generated by  elements of the form
$\delta f_i-\delta_Xf_i$ where $f_i\in \Cal{O}(X)$ generate $\Cal{O}(X)$ as $R$--algebra.}

\smallskip

We use the word ``flow'' in this context in order to suggest that in our main applications
we consider the $p$--adic axis as an arithmetic version of the time axis.

\smallskip

The derivation $\d_X$ is completely determined by its $\d$--flow. If $\d_X$ corresponds to the section
$s:\,J^0(X)\to J^1(X)$  of 
 $J^1(X)\to J^0(X)$ then
 $Z(\Cal{E}(\d_X))\subset J^1(X)$ coincides with the image of ths section.
 One easily checks that
if $\Cal{H}\in \Cal{O}(X)\,\widehat{}$ is such that $\d_X \Cal{H}=0$ then $\Cal{H}$ is a prime integral for $\Cal{E}(\d_X)$.
 All of the above can be transposed to the case when $X$ a $p$--formal scheme, locally a $p$--adic completion of a smooth scheme over $R$.

\smallskip

In what follows we choose  a smooth affine scheme $Y$ and  apply the constructions
discussed above to $X:=J^1(Y)$. In this case one can define a special class of $\d$-flows on $J^1(Y)$ which will be called {\it canonical $\d$-flows}.

\medskip

{\bf 1.3.3. Definition.} {\it A canonical $\d$--flow is  a $\d$--flow $\Cal{E}(\d_{J^1(Y)})$ with the property that the composition of $\d_{J^1(Y)}:\Cal{O}(J^1(Y))\to \Cal{O}(J^1(Y))$ 
with the pull back map $\Cal{O}(Y)\to \Cal{O}(J^1(Y))$ equals the universal
$p$--derivation $\d:\Cal{O}(Y)\to \Cal{O}(J^1(Y))$.}

\medskip

Notice  that in view of the universality property of $p$--jet spaces,  one gets a natural
closed embedding $\iota:\,J^2(Y)\to J^1(J^1(Y)).$
This induces  an injective  map (which we view as an identification) from the set of sections of $J^2(Y)\to J^1(Y)$ to the set of sections
of $J^1(J^1(Y))\to J^1(Y)$. The sections of $J^2(Y)\to J^1(Y)$ are in a natural bijection with canonical $\d$--flows on $J^1(Y)$ whereas  the sections
of $J^1(J^1(Y))\to J^1(Y)$ are in a bijection with (not necessarily canonical) $\d$--flows on $J^1(Y)$. 

\smallskip

Finally, consider   a system of arithmetic differential equations of order $2$,
$\Cal{F}\subset \Cal{O}(J^2(Y))$.

\medskip

{\bf 1.3.4. Definition.} {\it $\Cal{F}$ defines a $\d$--flow on $J^1(Y)$ if 
the map $Z^2(\Cal{F})\to J^1(Y)$ is an isomorphism.}

\medskip

In this  case then $Z^2(\Cal{F})\to J^1(Y)$ defines a section of $J^2(Y)\to J^1(Y)$ and hence a canonical $\d$--flow
$\Cal{E}(\d_{J^1(Y)})$ on $J^1(Y)$ such that 
$\iota(Z^2(\Cal{F}))=Z^1(\Cal{E}(\d_{J^1(Y)})).$

\smallskip
The differential algebra counterpart of the above definition yields the natural concept of flow on the (co)tangent space defining a second order differential equation.

\bigskip

\centerline{\bf 2.~ Painlev\'e VI and differential characters of elliptic curves}

\medskip

{\bf 2.1. Classical case.} The  family of  sixth Painlev\'e equations
depends on  four arbitrary constants $({\alpha ,\beta , \gamma ,\delta})$
and is classically written as
$$
\frac{d^2X}{dt^2}=\frac{1}{2}\left(
\frac{1}{X}+\frac{1}{X-1}+\frac{1}{X-t}\right)
\left(\frac{dX}{dt}\right)^2 -
\left(
\frac{1}{t}+\frac{1}{t-1}+\frac{1}{X-t}\right)\frac{dX}{dt} +
$$
$$
+\frac{X(X-1)(X-t)}{t^2(t-1)^2}
\left[ \alpha +
\beta\frac{t}{X^2}+\gamma\frac{t-1}{(X-1)^2}+
\delta\frac{t(t-1)}{(X-t)^2}\right]. \eqno{(2.1)}
$$
As R.~Fuchs remarked in 1907, (2.1) can be rewritten as the differenttial
equation for a (local) section $P:= (X(t),Y(t))$ of the generic elliptic curve
$E =E(t):\ Y^2=X(X-1)(X-t)$:
$$
t(1-t)\left[t(1-t)\frac{d^2}{dt^2}+(1-2t)\frac{d}{dt}-\frac{1}{4}\right]
\int_{\infty}^{(X,Y)}\frac{dx}{\sqrt{x(x-1)(x-t)}}=
$$
$$
=\alpha Y+\beta\frac{tY}{X^2}+\gamma\frac{(t-1)Y}{(X-1)^2}
+(\delta -\frac{1}{2})\frac{t(t-1)Y}{(X-t)^2}  \eqno{(2.2)}
$$
The l.~h.~s. of (2.2) can be called {\it the additive differential character $\mu$ of order two}
of $E$:  it is a non--linear differential expression in coordinates of $P$
such that $\mu (P+Q)=\mu (P) +\mu (Q)$ where $P+Q$ means addition
of points of the generic elliptic curve $E$, with infinity as zero. In particular,
$\mu (Q) =0$ for points of finite order. The point is that the integral in the l.h.s.
of (2.2) already has such additivity property, but it is defined only modulo periods,
and the latter are annihilated by the Gauss differential operator.

\smallskip

Thus $\mu (P)$  is defined up to multiplication by an invertible function of $t$.
If we choose a differential of the first kind $\omega$ on the generic curve and 
 the symbol of the Picard--Fuchs operator of the second order
annihilating periods of $\omega$,  the character will be defined uniquely.
In particular, if we pass to the analytic picture replacing the algebraic family
of curves $E(t)$ by the analytic  one $E_{\tau}:= \bold{C}/(\bold{Z}+\bold{Z}\tau )\mapsto \tau\in H$,
and denote by $z$ a fixed coordinate on $\bold{C}$, then (2.1) and (2.2) can be equivalently written
in the form
$$
\frac{d^2z}{d\tau^2}=\frac{1}{(2\pi i)^2}\sum_{j=0}^3
\alpha_j\wp_z(z+\frac{T_j}{2},\tau )
\eqno{(2.3)}
$$
where $(\alpha_0,\dots ,\alpha_3):=(\alpha ,-\beta ,\gamma ,\frac{1}{2}-\delta )$ and
$$
\wp (z,\tau ):= \frac{1}{z^2}+\sum_{(m,n)\ne (0,0)}
\left(\frac{1}{(z+m\tau +n)^2}-\frac{1}{(m\tau +n)^2}\right) .     
\eqno{(2.4)}
$$
Moreover, we have
$$
\wp_z (z,\tau )^2=4(\wp (z,\tau )-e_1(\tau ))(\wp (z,\tau)-e_2(\tau ))
(\wp (z,\tau )-e_3(\tau ))
\eqno{(2.5)}
$$
where
$$
e_i(\tau )=\wp (\frac{T_i}{2},\tau ),\ (T_0,\dots ,T_3)=
(0,1,\tau ,1+\tau )
 \eqno{(2.6)}
$$
so that $e_1+e_2+e_3=0$.

\smallskip

For a more geometric description, cf. sec. 4.5 below.

\medskip

{\bf  2.2. $p$--adic differential additive characters.}  The form (2.3) is more suggestive
than (2.1) for devising a $p$--adic version of PVI.
In purely algebraic terms, $z$ can be described as the logarithm of the
formal group law, and the rhs of (2.3) is simply a linear combination
of shifts of  $\wp_z$- (or $Y$-)coordinate by sections of the second order.

\smallskip

More precisely, using basic conventions of  [Bu05], let $p\ge 5$ be a prime,
$k$ an algebraic closure of $\bold{Z}/p\bold{Z}$, $R=W(k)$ the ring of $p$--typical
Witt vectors., as in 1.1.

\smallskip
Consider a  smooth projective  curve of genus one $E$ over $R$, with four
marked and numbered $R$--sections $P_i$, $i=0,\dots ,3,$  such that all divisors 
$2(P_i-P_j)$ are principal ones. Choosing, say, $P_0$ as zero section,
we may and will identify $E$ with its Jacobian 
and represent $E$ 
as the closure of the affine curve $y^2=4x^3+ax+b=4(x-e_1)(x-e_2)(x-e_3)$, with
$e_i\in R$, corresponding to $P_i-P_0$.
Assume that $E$  is {\it not a canonical lift
of its (good) reduction}, that is, does not admit a lift  of the Frobenius morphism of $E\otimes k$.

\smallskip
Choose the differential $\omega := dx/y$. Whenever we work with several elliptic curves
simultaneously, we will write $x_{E,\omega}, y_{E,\omega}$ in place of former $x,y$ etc.

\smallskip
The curve $E$ has a canonical $p$--differential character $\psi_{E,\omega}$  of order 2 ([Bu05], pp. 201 and 197),
which corresponds to $(2\pi i)^2 d^2z/d\tau^2$ in (2.3).

\medskip

{\bf 2.3. Painlev\'e $p$--adic equation and the problem of constants.} Now we can
directly write a $p$--adic version of (2.3) as
$$
\psi_{E,\omega} (Q)= \sum_{j=0}^3 \alpha_j s_j^*(y (Q))
\eqno(2.7)
$$
where $Q$ is a  variable section of $E/R$, and $s_j:\,E_j\to E_j$ is the shift by $P_j$.

\smallskip

At this point we have to mention two problems.

\smallskip

(A) In (2.3), $\alpha_j$ must be absolute constants rather than, say, functions of $t$.
Directly imitating this condition, we have to postulate that in (2.7), $\alpha_j$ must be
roots of unity or zero, i.~e. $\delta_p$--constants. It is desirable
to find a justification of such a requirement (or a version of it) in a more extended $p$--adic 
theory of Painlev\'e VI,  e.~g. tracing its source to the arithmetic analog
of isomonodromy deformations.

\smallskip

(B) The relevance and non--triviality of the problem of  ``constants'' in the
$p$--adic differential equations context is also implicit in our exclusion of
 those $E$ that are canonical lifts of their reductions.

\smallskip
A formal reason for this exclusion was the fact for such $E$ the basic differential character $\psi$ has order 1 rather than 2,
thus being outside the framework of Painlev\'e VI. But the analogy with functional case
suggests that canonical lifts should be morally considered
as analogs of families with constant $j$--invariant in the functional case.

\smallskip

This agrees also with J.~Borger's philosophy that   Frobenius lift(s)
should be treated as descent data to an {\it  ``algebraic geometry below
$\roman{Spec}\,\bold{Z}$''}  (cf. [Bo05]). In our particular case 
the latter might be called ``geometry over $p$--typical field with one element''
(Borger's suggestion in e-mail to Yu.M. of April 22, 2013).

\smallskip

Indeed, canonical lifts $X$ in such a geometry are endowed with an isomorphism 
$X^{\phi}\to X^{(p)}$ that can be seen as a categorification
of the identity  $c^{\phi}= c^p$ defining roots of unity.
\smallskip

In the $p$--adic case, however, if we decide to declare $j$--invariants
of canonical lifts ``constants'' in some sense, this will require 
a revision of the latter notion. These invariants generally are not
roots of unity: cf. Finotti's papers in http://www.math.utk.edu/ $\widetilde{}$ finotti/ .

\bigskip

\centerline{\bf 3. Symmetries and variants}

\medskip

{\bf 3.1. Lemma.} {\it (i) (Landin's transform). In the notations of 2.3, denote for each $i=1,2,3$ by 
$\pi_i:\,E\to E_i:=E/\langle P_i\rangle$ the respective isogeny. Let $\omega_i$
be the 1--form on $E_i$ such that $\pi^*(\omega_i)=\omega$. Then
$$
\pi_i^*(y_{E_i,\omega_i})=y_{E,\omega}+s_i^*(y_{E,\omega}).
\eqno(3.1)
$$ 

(ii) We have
$$
\pi_i^*(\psi_{E_i,\omega_i}) =  \psi_{E,\omega}.
\eqno(3.2)
$$

}
{\bf Proof.} (i) If we choose an embedding of $R$ in $\bold{C}$ identifying
$E(\bold{C})$ with $\bold{C}/(\bold{Z}+\bold{Z}\tau )$ in such a way that $P_i$ becomes
the point $\tau/2$ (modulo periods), then (3.1) turns into the classical Landin identity (cf. [Ma96], sec. 1.6):
$$
\wp (z,\tau/2 )  = \wp (z,\tau )   +  \wp (z+\tau/2,\tau )
$$

(ii) The identity (3.2) can be stated and proved  in  wider generality. Namely, let $E,E^{\prime}$ be
two smooth elliptic curves over $R$, not admitting lifts of Frobenius.
Let $\pi :\,E\to E^{\prime}$ be an isogeny of degree prime to $p$ and
and $\omega, \omega^{\prime}$ such bases of 1--forms on $E,E^{\prime}$ that
$\pi^*(\omega^{\prime})=\omega$. Then 
$$
\pi^*(\psi_{E^{\prime},\omega^{\prime}}) =  \psi_{E,\omega}.
\eqno(3.3)
$$
In fact from [Bu05], Theorem 7.34, it follows that 
$\pi^*(\psi_{E^{\prime},\omega^{\prime}}) = c \psi_{E,\omega}$ for some $c\in R$.

\smallskip
Now  consider the second jet space $J^2(E)$. Then $\psi =\psi_{E,\omega}$ is an element 
of $\Cal{O}(J^2(E))$ (cf. [Bu05], p. 201). Denote by $\phi^i:\, J^i(E)\to J^0(E)$ the map defined
at the end of sec. 1.2. Put 
$$
\omega^{(i)}:= \frac{(\phi^i)^{*}(\omega )}{p^i},\ i=0,1,2.
\eqno(3.4)
$$
Then in view of [Bu05], p. 203 (where our $\omega^{(i)}$ were denoted $\omega_i$), we have
$$
d\psi =p\omega^{(2)} +\lambda_1\omega^{(1)}+\lambda_0\omega^{(0)}
\eqno(3.5)
$$
with $\lambda_1 \in R$, $\lambda_1 \in R^{\times}$ depending on $(E,\omega )$.  Notice that here $d$ means 
the usual differential taken in the ``geometric'', or ``vertical'' direction, that is $dc=0$ for any $c\in R$.
Moreover, in [Bu97] it was proved that if $E$ is defined over $\bold{Z}_p\subset R$, 
then $\lambda_0=1$, $\lambda_1=-a_p$ where
$a_p$ is the trace of Frobenius on the reduction $E\,\roman{mod}\,p$.

\smallskip

Returning to the proof of (3.3), we can now  compare the $\omega^{(2)}$--contributions 
to $d \psi_{E^{\prime},\omega^{\prime}}$ in two different ways. First, we have
$$
\pi^*(d \psi_{E^{\prime},\omega^{\prime}})=\pi^*(p(\omega^{\prime})^{(2)}+\dots )= p\,(\pi^*(\omega^{\prime})^{(2)}) + \dots
$$
Second,
$$
\pi^*(d \psi_{E^{\prime},\omega^{\prime}})= d \pi^*( \psi_{E^{\prime},\omega^{\prime}})=
d(c\,\psi_{E,\omega})=c\,\,(p\omega^{(2)}+\dots ).
$$
This shows that $c=1$.

\medskip

{\bf 3.2. Two versions of PVI.}  Put $Y:= E\setminus \cup_{i=0}^3P_i.$
Denote by $r:= \sum_{j=0}^3 \alpha_j s_j^*(y)\in \Cal{O}(Y)$ (cf. (2.7)).  

\smallskip
Below we will denote by $\rho$ {\it either $r$, or $\phi( r ) \in\Cal{O}(J^1(Y))$.}
Although the introduction of the version with $\phi ( r )$ is not motivated
at this point, we will see in sec. 4 below that exactly this version admits
a ``$p$--adic Hamiltonian'' description.
\smallskip
The character $\psi$ induces 
the map of sets, that we will also denote $\psi :\, E( R )\to R$; similarly, $\rho$ induces
the map of sets $\rho\,:\, Y( R )\to R$.

\medskip
{\bf 3.2.1. Proposition.} {\it Denote by the symbol 
$
PVI(E,\omega, P_0,P_1,P_2,P_3,\alpha_0,\alpha_1,\alpha_2,\alpha_3)
$
any of the two equations
$$
\psi ( P )=\rho ( P ).
\eqno(3.6)
$$

Let $\pi := \pi_2:\ E\to E^{\prime}:= E/\langle P_2\rangle$ as in Lemma 3.1.
Put $P_0^{\prime}=\pi (P_0)$ (zero point), $P_1^{\prime}=\pi (P_1)$, and 
choose remaining two points of order two $P_2^{\prime}, P_3^{\prime}$. Assume
that $\pi^*(\omega^{\prime})=\omega$.

In this case for any solution $Q$ to  (3.6), the point $Q^{\prime}:=\pi (Q)$
will be a solution to
$$
PVI^{\prime}=PVI(E^{\prime},\omega^{\prime}, P_0^{\prime},P_1^{\prime},P_2^{\prime},P_3^{\prime},\alpha_0,\alpha_1,
0,0),
\eqno(3.7)
$$ 
and conversely, if $Q^{\prime}$ is a solution to $PVI^{\prime}$, then $Q$ is a solution to $PVI$.
}
\medskip

{\bf Proof.} Let $\rho = \phi^j( r )$, $j=0$ or $1$. Our statement results from the following
calculation, using (3.1), (3.2):
$$
\psi_{ E^{\prime},\omega^{\prime}} (\pi (Q))=\psi_{E,\omega}(Q)=
$$
$$
\phi^j(\alpha_0(y_{E,\omega}(Q)+y_{E,\omega}(Q+P_2))
+ \alpha_1(y_{E,\omega}(Q+P_1)+y_{E,\omega}(Q+P_1+P_2)))=
$$
$$
\phi^j(\alpha_0y_{ E^{\prime},\omega^{\prime}} (\pi (Q))
+ \alpha_1y_{ E^{\prime},\omega^{\prime}} (\pi (Q)+\pi(P_1))).
$$

\medskip

Of course, similar statement will hold if $\pi_2$ is replaced by $\pi_1$ or $\pi_3$.
\medskip

{\bf  3.3. Two more versions of PVI.} The character $\psi$ is of course  not algebraic,
but the equation (3.5) shows that it has an algebraic vertical differential.
The same is obviously true for the rhs of (2.7): one easily sees that
if $E$ is given in the Weierstrass form $y^2=f(x)$, then
$$
dr = -\frac{1}{2}
\sum_{j=0}^3 \alpha_j s_j^*(f^{\prime}(x)))\, \omega.
\eqno(3.8)
$$
Hence we may consider two more versions of the arithmetic $PVI$: the condition of vanishing
of one of the following 1--forms on $J^2(E)$:
$$
d\psi-dr=p\omega^{(2)}+ \lambda_1\omega^{(1)}+ \left(\lambda_0-\frac{dr}{\omega}\right)\omega^{(0)},
\eqno(3.9a)
$$
$$
d\psi-d\phi ( r )=p\omega^{(2)}+\left(\lambda_1-p\phi\left(\frac{dr}{\omega}\right)\right)\omega^{(0)}+ \lambda_0\omega^{(0)}.
\eqno(3.9b)
$$

But in the $p$--adic situation  solutions to (3.6) and (3.9a,b) respectively
are not related to each other in the way we would expect by analogy with
usual calculus.
\smallskip
Indeed, let us consider a smooth function $f(x,y,y',...,y^{(r)})$ in $r+2$ variables defined on ${\bold{R}^{r+2}}$, where the latter is viewed as the $r$--th jet space of the first projection $\bold{R}^2\to \bold{R}$,
$(x,y)\mapsto x$. Let $u=u(x)$ be an unknown smooth function $u: \bold{R}\to \bold{R}$.
Then the equation 
$$
f(x,u(x),u'(x),...,u^{(r)}(x))=0
\eqno(3.10)
$$
is related to the $1$-form  
$df$ on $\bold{R}^{r+2}$
as follows. If $u$ solves (3.10), then taking the derivative of (3.10) with respect to $x$ we get
$$
\nabla^r(u)^*(df)=0
\eqno(3.11)
$$
where $\nabla^r(u): \bold{R}\to \bold{R}^{r+2}$ is given by $\nabla^r(u)(x)=(x,u(x),...,u^{(r)}(x))$.
\smallskip
However in the case of arithmetic jet spaces the situation is different.
\smallskip
 Indeed, let $f\in \Cal{O}(J^r(\bold{A}^1))=R[y,y',...,y^{(r)}]\widehat{\ }$ where $J^r(\bold{A}^1)$ is the arithmetic jet space of the affine line over $R$. Furthermore,   let $u\in R$ be a solution of
$$
f(u,\d u,...,\d^r u)=0
\eqno(3.12)
$$
where $\delta = \delta_p:R\to R$ is the standard $p$--derivation on $R$.  Applying $\delta$ to (3.12)
we get:
$$
\frac{1}{p}(f^{(\phi)}(u^p+p\d u, (\d u)^p+p\d^2 u,...,(\d^r u)^p+p\d^{r+1} u)-f(u,\d u,...,\d^r u)^p)=0
\eqno(3.13)
$$
where  $f^{(\phi)}$ is $f$ with coefficients twisted by the Frobenius  $\phi$. To compare
this latter equation with (3.11), one can apply
 the Taylor formula and get
$$
 \frac{f^{(\phi)}(\nabla^r(u)^p)-f(\nabla^r(u))^p}{p}
 + \sum_{i=0}^r\frac{\partial f^{(\phi)}}{\partial y^{(i)}}(\nabla^r(u)^p)(\d^{i+1} u)+M=0,
 \eqno(3.14)
 $$
where $\nabla^r(u)=(u,\d u,...,\d^r u)$.

\smallskip
The first term  here is an arithmetic version of the pullback of the term $\frac{\partial f}{\partial x}dx$ in $df$.   The second term clearly  involves $df$. However   $M$ involves higher partial derivatives of $f$; it is highly non-linear, it is divisible by $p$, and ``the more non-linear its terms are the more they are divisible by $p$".  In some sense, since $p$ is small, one can view (3.14) as a non-linear deformation of (3.11).  
\medskip
We will continue discussion of this formalism in sec. 5 below. The reader may wish to
skip the next section, or to postpone reading it.

\bigskip

\centerline{\bf 4. Hamiltonian formalism}

\medskip

{\bf 4.1. Classical case.} The classical PVI equation written in the form (2.3)
can be represented as a Hamiltonian flow on the variable two--dimensional phase space
(twisted cotangent spaces to a versal family of elliptic curves, cf. [Ma96] and the end of this section), with time--dependent Hamiltonian:
$$
\frac{dz}{d\tau}=\frac{\partial\Cal{H}}{\partial y},\
\frac{dy}{d\tau}=-\frac{\partial\Cal{H}}{\partial z},
\eqno{(4.1)}
$$
where
$$
\Cal{H}:=\Cal{H} (\alpha_0,\dots ,\alpha_3):=\frac{y^2}{2}-\frac{1}{(2\pi i)^2}
\sum_{j=0}^{3}\alpha_j\wp (z+\frac{T_j}{2},\tau ).
\eqno{(4.2)}
$$

In more geometric terms, this means that  solutions to the PVI become leaves of the 
null--foliation of the following closed two--form:
$$
\omega =\omega(\alpha_0,\dots ,\alpha_3):=
2\pi i(dy\wedge dz-d\Cal{H}\wedge d\tau )=
$$
$$
=2\pi i(dy\wedge dz-ydy\wedge d\tau ) +
\frac{1}{2\pi i}
\sum_{j=0}^{3}\alpha_j\wp_z (z+\frac{T_j}{2},\tau )dz\wedge d\tau .
\eqno{(4.3)}
$$
The extra factor $2\pi i$ makes $\omega$  defined over $\bold{Q}[\alpha_i]$ on a natural
algebraic model of (twisted) relative cotangent bundle to the respective versal
family of elliptic curves.
\smallskip

In the expression (4.3), the summand $2\pi i\,dy\wedge dz$ is the canonical fibrewise symplectic
form on the relative cotangent bundle. The terms involving $d\tau$ uniquely determine the differential
of the (tme dependent)  Hamiltonian.

\smallskip

Moreover, $\omega$ is not just closed, but is a global differential: $\omega = d\nu$
where the form
$$
\nu =\nu (\alpha_0,\dots ,\alpha_3) :=2\pi i\,(ydz-\frac{1}{2}y^2d\tau)
+d\roman{log}\,\theta(z,\tau )
+2\pi i\,G_2(\tau )d\tau +
$$
$$
+\frac{1}{2\pi i}\sum_{j=0}^3\alpha_j\wp (z+\frac{T_j}{2},\tau )d\tau
\eqno{(4.4)}
$$
also descends to an appropriate algebraic model, and the Hamiltonian $\Cal{H}$ is again encoded
in the $d\tau$--part of $\nu$. Here is a convenient way to represent this encoding:
$$
\Cal{H}(\alpha_0,\dots \alpha_3) =i_{{}_{\partial_\tau}}\left(\frac{y^2}{2}d\tau+ \frac{1}{2\pi i} (\nu (0,\dots ,0)-
\nu (\alpha_0,\dots ,\alpha_3))\right),
\eqno(4.5)
$$
where  $\partial_\tau := \frac{\partial}{\partial \tau}$.

\smallskip

Finally, the last summand in (4.3) is simply the additive differential character $2\pi i \dfrac{d^2z}{d\tau^2}$
that is generally denoted $\psi$ in the $p$--adic case.

\smallskip

In the following subsections, we try to imitate this description of Hamiltonian structure
for $p$--adic PVI equations.
The reader should be aware that our treatment is somewhat {\it ad hoc},
and must be considered as a tentative step towards a more coherent vision
of Hamiltonian flows with $p$--adic time. In fact, we do  not yet
have an appropriate version of $dp$ replacing $d\tau$ and generally
do not know what are differential forms involving ``differentials in the arithmetical direction''.

\medskip
{\bf 4.2. Arithmetical case: preparation.} Let  $Y$ be an formal affine scheme over $R=W(k)$. Modules of vertical differential
forms on  $Y$ are defined as
$$
\Omega_Y = \roman{lim\,inv}\,\Omega_{Y_n/R_n}
$$
where $R_n=R/p^{n+1}R$, $Y_n=Y\otimes_R R_n$.

\smallskip

Let now $Z\subset J^n(Y)$ be a closed formal subscheme defined by the ideal $I_Z\subset \Cal{O}(J^n(Y)). $ 
Put
$$
\Omega^{\prime}_Z:= \frac{\Omega_{J^n(Y)}}{\langle I_Z \Omega_{J^n(Y)}, dI_Z\rangle}
\eqno(4.6)
$$

Given a system of arithmetic differential equations $\Cal{F}\subset \Cal{O}(J^r(Y))$, denote by $Z^r:=Z^r(\Cal{F})$
the ideal generated by $\Cal{F}$. For each $s\le r$, there is a natural map
$\pi_{r,s} :\, Z^r\to J^s(Y).$

\smallskip

Generally, the natural maps $\phi^*$ respect degrees of differential forms, 
one can define natural maps 
$
\phi^*/p^i :\, \Omega^i_{J^{r-1}(Y)}\to \Omega^i_{J^{r}(Y)}
$
and, for $f\in \Cal{O}(J^2(Y))$, they induce maps which we will denote
$$
\frac{\phi^*_Z}{p^i} :\, \Omega^i_{J^1(Y)}\to \Omega^{\prime i}_{Z^2(f)}
\eqno(4.7)
$$

\medskip

{\bf 4.2.1. Definition.} {\it   We say that   $\Cal{F}\subset \Cal{O}(J^r(Y))$   defines a generalised canonical
$\delta$--flow on $J^s(Y)$, if the
induced map
$$
\pi_{r,s}^*\Omega_{J^s(Y)}\to \Omega^{\prime}_{Z^r}
$$
is injective, and its cokernel is annihilated by a power of $p$.}

\medskip

The cokernel here intuitively measures ``how singular'' $\Cal{F}$ is on the closed fibre of $Y$.

\medskip
{\bf 4.2.2. Definition.} {\it a) Let $X$ be a smooth surface over $R$ (or the $p$--adic completion of such a surface). 

A symplectic form on $X$ is an invertible $2$--form on $X$. 

A contact form on $X$ is an $1$--form on $X$ such that $d\nu$ is symplectic.

\smallskip

b) Let $Y$ be a smooth curve over $R$. An $1$--form on $X:=J^1(Y)$ is called canonical,
if $\nu=f\beta$, where $f\in \Cal{O}(X)$ and $\beta$ is an $1$--form lifted
from $Y$.}
\medskip

Notice that any closed canonical $1$--form on $X=J^1(Y)$ is lifted from $Y$.

\medskip

We now come to the main definition.

\medskip

{\bf 4.3. Definition.} {\it  Let  $Y$ be a smooth affine curve over $R$ and let $f\in \Cal{O}(J^2(Y))$ 
be a function defining  a generalised canonical $\d$--flow on $J^1(Y)$. 

\smallskip
a) The respective generalized $\d$--flow is called  Hamiltonian with respect to the symplectic form $\eta$ on $J^1(Y)$,
 if $\phi^*_Z\eta =\mu\cdot \eta$
in $\Omega^{\prime 2}_{Z^2(f)}$ for some $\mu \in pR$ called the eigenvalue.

\smallskip

b) Assume that moreover $\eta =d\nu$ for some canonical $1$--form  $\nu$ on $J^1(Y)$.
Then we call
$$
\epsilon:= \frac{\phi_Z\nu-\mu \nu}{p}\in \Omega^{\prime}_{Z^2(f)}
\eqno(4.8)
$$
an Euler--Lagrange form.}

\medskip

We consider (4.8) as an (admittedly, half--baked) arithmetical analog of the expression $i_{{}_{\partial_\tau}}(ydz-\Cal{H} d\tau )$
(cf. (4.5))
in the same sense as the $p$--derivation
$$
\d_p (x)=\frac{\phi (x)-x^p}{p}
$$
is an analog of $\partial_{\tau}$.

\medskip

Now we pass to the arithmetical PVI. Let again $E$ be an elliptic curve over $R$ that does not admit a lift of Frobenius and let $\psi\in \cO(J^2(E))$ be the canonical $\d$-character of order $2$ attached to an invertible 
$1$--form $\omega$ on $E$. Consider the symplectic form  $\eta=\omega^{(0)}\wedge \omega^{(1)}$ on $J^1(E)$:
cf. (3.4). Let $Y\subset E$ be an affine open set and let $r\in \cO(Y)$. Assume in addition that $Y$ has an \'{e}tale coordinate. (The
basic example is $E$ with sections of the second order deleted). 

\smallskip

Denoting such an \'etale coordinate by $T$, put $\Cal{A}_2 = K[[T,T^{\prime}]]$,
$\Cal{A}_3 = K[[T,T^{\prime},  T^{\prime\prime}]]$, where $K$ is the quotient ring of $R$.

\medskip

{\bf  4.4. Proposition.} {\it 
The following assertions hold:
\smallskip
1) The function  $f=\psi-\phi(r)$  defines a generalised canonical $\d$--flow on $J^1(Y)$ which  is Hamiltonian with respect to $\eta$.
\smallskip
2)  There exists a canonical  $1$-form $\nu$ on $X$ such that $d\nu=\eta$; in particular the symplectic form $\eta$ is exact and if  $\epsilon:=\frac{1}{p}(\phi_Z^*\nu-\mu\nu)$ is the Euler-Lagrange form then $p\epsilon$ is closed.
 \smallskip
 3) Let 
$\epsilon$ be Euler-Lagrange form and $f=\psi-\phi(r)$. Then we have the following equality in $\Omega_{\Cal{A}_3}$:
$$
\epsilon=f\omega^{(1)}-\frac{1}{p}(\phi^*-\mu)\nu.
$$
\smallskip
4) Let $r_1,r_2\in \cO(Y)$ be such that
 $r_2-r_1=\partial s$, for some $s\in \cO(Y)$,
 where $\partial$ is the derivation on $E$ dual to $\omega$. (This holds for  two right hand sides of any two  PVI equations). Consider the equations  $\psi-\phi(r_1)$ and $\psi-\phi(r_2)$ respectively. Then there exists a canonical $1$--form $\nu$ on $J^1(Y)$ such that $d\nu=\eta$ and such that, if $\epsilon_1,\epsilon_2$ are the corresponding Euler--Lagrange forms, then  
 :}
  $$
  \epsilon_1-\epsilon_2=\frac{1}{p}d\phi(s)  \in \Omega_{\Cal{A}_2}.   
  $$  
\medskip

{\bf Remark.} We will deduce from this Proposition below (cf. Corollary 5.3.3) that in fact $p$--adic PVI
in the form (3.6) defines a generalised canonical $\delta$--flow. However, it does not define a $\delta$--flow in the 
sense of  Definition 1.3.4. This  motivated our Definition 4.2.1.

\medskip
{\it Proof}.   From (3.4), 
we get the following equality in $\Omega^2_{J^2(Y)}$:
$$
\frac{\phi^*\eta}{p^2} =\frac{\phi^*\omega^{(0)}}{p}\wedge \frac{\phi^*\omega^{(1)}}{p}=\omega^{(1)}\wedge \omega^{(2)}.
$$
Recall the formula (3.9b): 
$$
d f = p\omega^{(2)}+(\lambda_1-p\phi (\partial r))\omega^{(1)}+\lambda_0\omega^{(0)}
$$
 in $\Omega_{J^2(E)}$ where $\lambda_1\in R$, $\lambda_0\in R^{\times}$. Hence, if we keep notation $\omega^{(0)},\omega^{(1)}$ also for the images of the respective forms in $\Omega^{\prime}_{Z^2(f)}$, in view of
$\omega^{(1)}\wedge df=0$ in $\Omega^{\prime}_{Z^2(f)}$, 
we have the following equality in $\Omega^{\prime 2}_{Z^2(f)}$:
$$
\phi^*_Z \eta  =  -p\cdot \omega^{(1)}\wedge ((\lambda_1-p\phi(dr))\omega^{(1)}+\lambda_0\omega^{(0)})= p\lambda_0 \eta,
$$
This completes the proof of 1).

\smallskip

 Write now $\omega=d L=\frac{dL}{dT}dT$ where $L=L(T)\in TK[[T]]$. (For instance, we can take $L$ to be the formal logarithm of $E$).
Then
$$
{\Cal{A}_3}=K[[T,\phi(T),\phi^2(T)]]=K[[L,\phi(L),\phi^2(L)]].
$$
So the image of $\psi$ in ${\Cal{A}_3}$ is
$$
\psi=\frac{1}{p}(\phi^2(L)+\lambda_1\phi(L)+p\lambda_0L)+\lambda_{-1}
$$
for some $\lambda_{-1}\in R$.
Therefore the maps

$$
\Omega_{{\mathcal A}_2}\ra \frac{\Omega_{{\mathcal A}_3}}{\langle f\Omega_{{\mathcal A}_3},df\rangle},\ \ \ \ \Omega^2_{{\mathcal A}_2}\ra \frac{\Omega^2_{{\mathcal A}_3}}{\langle f\Omega^2_{{\mathcal A}_3},df\wedge \Omega_{{\mathcal A}_3}\rangle}
$$
are isomorphisms and so we have induced Frobenii maps $\phi_f^*:\Omega_{{\mathcal A}_2}\ra \Omega_{{\mathcal A}_2}$ and $\phi_f^*:\Omega^2_{{\mathcal A}_2}\ra \Omega^2_{{\mathcal A}_2}$.
 Since $T$ is \'{e}tale the lift of Frobenius $T\mapsto T^p$ on $\bold{A}^1$ extends to a lift of Frobenius $\phi_0:\widehat{Y}\ra \widehat{Y}$ of the $p$--adic completion of $Y$. Also the derivation 
 $\frac{d}{dT}$ on $R[T]$ extends to a derivation still denoted by $\frac{d}{dT}$ on $\cO(\widehat{Y})$. We claim that
$$\frac{\phi(L)-\phi_0(L)}{p},$$
which a priori is an element of ${\mathcal A}_2$, actually belongs to $\cO(J^1(Y))$. Indeed we have the following expansion in ${\mathcal A}_2$:
$$
\frac{\phi(L)-\phi_0(L)}{p}  =  \frac{L^{(\phi)}(T^p+pT')-L^{(\phi)}(T^p)}{p}
$$
$$
= \sum_{i=1}^{\infty} \frac{p^{i-1}}{i!}      \frac{d^i L^{(\phi)}}{dT^i}(T^p)(T')^i
 =\sum_{i=1}^{\infty} \frac{p^{i-1}}{i!}\phi_0\left( \frac{d^iL}{dT^i}\right) (T')^i
$$
$$
= \sum_{i=1}^{\infty} \frac{p^{i-1}}{i!}\phi_0\left( \left(\frac{d}{dT}\right)^{i-1}\left(\frac{\omega}{dT}\right)\right) (T')^i \in
 \cO(J^1(Y)),
 $$
where the superscript $(\phi)$ means twisting coefficients by $\phi$. The latter inclusion follows because $T'=\d T\in \cO(J^1(Y))$,  $\frac{\omega}{dT}\in \cO(Y)\subset \cO(\widehat{Y})$, and the latter is stable under $\frac{d}{dT}$ and $\phi_0$.
Now set
$$\nu:=-\frac{\phi(L)-\phi_0(L)}{p}\,\omega\in \cO(J^1(Y))\,\omega\in \Omega_{J^1(Y)}.$$
Then
$$
d\nu  =  -d\left(\frac{\phi(L)-\phi_0(L)}{p}\right)\wedge \omega
$$
$$
= -d\left(\frac{\phi(L)}{p}\right)\wedge \omega+d\left(\frac{\phi_0(L)}{p}\right)\wedge \omega
$$
$$
= -\omega^{(1)}\wedge \omega^{(0)}= \eta .
$$
This completes the proof of 2).  
\smallskip
Next for $f=\psi-\phi(r)$ we have the following computation in $\Omega_{{\mathcal A}_3}$:
$$
p \epsilon  =  \phi_f^*\nu-\mu\nu
$$
$$
= -\phi_f^*\left( \frac{\phi(L)-\phi_0(L)}{p} \omega\right) +\mu\frac{\phi(L)-\phi_0(L)}{p}\omega
$$
$$
= -\phi_f \phi(L)\omega^{(1)}+\phi\phi_0(L)\omega^{(1)}+
\mu\frac{\phi(L)-\phi_0(L)}{p}\omega
$$
$$
=(\lambda_1\phi(L)+p\lambda_0L+\phi\phi_0(L)-p\phi(r)+p\lambda_{-1})\omega^{(1)}+\mu\frac{\phi(L)-\phi_0(L)}{p}\omega^{(0)}
$$
$$
=(\lambda_1\phi(L)+p\lambda_0L+\phi\phi_0(L)-p\phi(r)+p\lambda_{-1})\omega^{(1)}
$$
$$
 + \phi^*\left(\frac{\phi(L)-\phi_0(L)}{p}\omega^{(0)}\right)-(\phi^*-\mu)\nu
 $$
 $$
=\left(\phi^2(L)+\lambda_1\phi(L)+p\lambda_0L-p\phi(r)+p\lambda_{-1}
\right)\omega^{(1)}-(\phi^*-\mu)\nu
$$
$$
=pf\omega^{(1)}-(\phi^*-\mu)\nu.
 $$
This ends the proof of assertion 3). Assertion 4) follows from the fact that
$$\epsilon_1-\epsilon_2=\phi(\partial s)\omega^{(1)}=\frac{1}{p}d(\phi(s)).$$

\medskip

{\it Remarks. a)} Some of our arguments above break down if $\psi-\phi(r)$ is replaced by $\psi-r$. This is our main motivation for studying $\psi-\phi(r)$.
\smallskip
{\it b)} Assertion 3) implies that if
$[\epsilon]_2,[f\omega^{(1)}]_2\in \frac{\Omega_{{\mathcal A}_3}}{(\phi^*-\mu)\Omega_{{\mathcal A}_2}}$ are the  images 
of $\epsilon,f\omega^{(1)}$ then
$$
[\epsilon]_2=[f\omega^{(1)}]_2.
$$
This justifies our suggestion that $f$ is the ``Euler--Lagrange equation attached to our Hamiltonian data". Assertion 4) says that the Euler--Lagrange forms of various PVI equations differ by exact forms.
\smallskip

{\it c)} Finally, we could treat in this way also the multicomponent version of PVI and the degeneration PV
as they are described in [Ta].

\bigskip

{\bf 4.5. Geometry of PVI in various categories.}
Below we essentially reproduce from [Ma96] a geometric description of the natural habitats of various forms of PVI
including its arithmetic version.

\smallskip

Our series of constructions starts with
a non--constant family (``pencil'') of elliptic curves over one--dimensional base.
in one of several natural categories: schemes over a field of characteristic zero
or a ring of algebraic integers, or a $p$--adic completen of the latter,
analytic spaces etc.
\smallskip

a. Let $(\pi :E\to B;D_0,\dots ,D_3)$ be a pencil of
compact smooth curves of genus one, with variable absolute invariant,
endowed with four labelled sections $D_i$ such that
if any one of them is taken as zero, the others
will be of order two.

\smallskip

We will call $E$ {\it a configuration space} of PVI
(common for all values of parameters.) Solutions to
all equations will be represented by some multisections
of $\pi .$

\smallskip

b. Let $\Cal{F}$ be the subsheaf of the sheaf of
vertical 1-forms $\Omega^1_{E/B}(D_3)$ on $E$
with pole at $D_3$ and residue 1 at this pole.
It is an affine twisted version of $\Omega^1_{E/B}$
which is the sheaf of sections of the relative
cotangent bundle $T^*_{E/B}.$ Similarly,
$\Cal{F}$ itself ``is'' the sheaf of sections
of an affine line bundle $F=F_{E/B}$ on $E.$
More precisely, we can construct  a bundle
$\lambda :F\to E$ and a form $\nu_F\in\Gamma (F, \Omega^1_{F/B}(\lambda^{-1}(D_3))$
such that the map
$$
\{\roman{local\ section}\ s\ \roman{of}\ F\}\mapsto
s^*(\nu_F)
$$
identifies the sheaf of sections of $F/E$ with $\Cal{F}.$

\smallskip

We will call $F$ {\it a phase space} for PVI (again, common
for all parameter values.) It is this space that carries a canonical symplectic
form (relative over the base) rather than  the usual cotangent bundle.

\smallskip

c. $E$ carries a distinguished family of algebraic/arithmetic curves
transversal to the fibers of $E$: considered as multisections of
$E/B$ they are of finite order (if any of $D_i$ is chosen as zero.)
It is important that each curve of this family has a canonical
lifting to $F$ (for its description, see [Ma96], especially
(2.12) and (2.29).)

\smallskip

d. $F$ carries a closed 2-form $\omega$ which can be characterised
by the following two properties:

\smallskip

{\it i). The vertical part of $\omega$, i.~e. its restriction
to the fibers of $\pi\circ\lambda : F\to B$, coincides
with $d_{F/B}(\nu_F).$}

\smallskip

{\it ii). Any canonical lift to $F$ of a connected multisection
of finite order
of $E\to B$, referred to above, is a leaf of the
null--foliation of $\omega$.}

\smallskip

e. $E$ also carries four distinguished closed two--forms
$\omega_0,\dots ,\omega_3.$ They are determined,
up to multiplication by a constant, by the following
properties.

\smallskip

{\it iii). The divisor of $\omega_i$ is $\dfrac{D_jD_kD_l}{D_i^3}$
where $\{ i,j,k,l\}=\{ 0,1,2,3\}.$}

\smallskip

{\it iv). Identify the sheaves $\Omega^2_E$ and $\pi^*(\Omega^1_{E/B})^{\otimes 3}$
on $E$ using the Kodaira--Spencer isomorphism
$\pi^*(\Omega^1_B)\cong (\Omega^1_{E/B})^{\otimes 2}$
and the exact sequence $0\to \pi^*(\Omega^1_B)\to
\Omega^1_E \to \Omega^1_{E/B} \to 0.$  Then the image of
$\omega_i$ in $\pi^*(\Omega^1_{E/B})^{\otimes 3}$
considered in the formal neighborhood of $D_i$ is the cube
of a vertical 1--form with a constant residue along $D_i.$}

\smallskip

Notice that up to introduction of $\omega$, all constructions were valid both in geometric
and arithmetic cases. It is precisely the absence of the differential
in arithmetic direction was the object of our concerns in 4.2--4.4 above.

\smallskip

The affine space $P_0:= \omega+\sum_{i=0}^3 \bold{C}\lambda^*(\omega_i)$
of closed two--forms on $F$ is our version of the
{\it moduli space of the PVI equations} replacing the
classical $(\alpha ,\beta ,\gamma ,\delta )$--space.

\smallskip

We can now summarize our geometric definition of PVI equations and
their solutions.

\smallskip

{\bf 4.5.1.  Definition} {\it a) A Painlev\'e two--form
on $F$ is a point $\Omega\in P_0.$
\smallskip
b) The Painlev\'e foliation corresponding to $\Omega$
is the null--foliation of $\Omega$.
\smallskip
c) The solutions to the respective Painlev\'e equation
are the leaves of this foliation (in the Hamiltonian
description). 
\smallskip
The form $\omega$ corresponds to
$(\alpha_0 ,\alpha_1,\alpha_2,  \alpha_3 )=(0,0,0,0).$}

\bigskip

\centerline{\bf 5. Pfaffian congruences}

\medskip

{\bf 5.1. Notation.} We will continue here the discussion started in 3.3.

For $u\in R$, define the following infinite vectors in $R\times R\times \dots $:
$$
\nabla(u) : = (u,\d u,\d^2 u,...),
$$
$$
 \d \nabla(u) :=  (\d u,\d^2 u,\d^3 u,...),
 $$
 $$
 \nabla(u)^p := (u^p,(\d u)^p,(\d^2 u)^p,...).
 $$
 
For $g\in R[y,y',...,y^{( r )}]\,\widehat{\ }$, denote by $g^{(\phi)}$ be the series obtained from $g$ by applying $\phi$ to the coefficients of $g$. Moreover, put
$$
dg :=  \sum_{i=0}^r \frac{\partial g}{\partial y^{(i)}} dy^{(i)},
$$
$$
\frac{\partial g}{\partial p}  :=  \frac{1}{p}(g^{(\phi)}(y^p, \dots (y^{(r)})^p)-g(y,\dots,y^{(r)})^p),
$$
$$
\frac{\partial\ \ \ }{\partial \nabla(y)}  :=   \left(\frac{\partial}{\partial y}, \frac{\partial}{\partial y'} ,\frac{\partial}{\partial y''}
\dots \right),
$$
Denoting by $\langle \, , \rangle$ the basic pairing between 1--forms and vector fields,
and writing   $\left( \frac{\partial}{\partial \nabla(y)}\right)^t$ for the  column  transpose of $ \frac{\partial}{\partial \nabla(y)}$,
we get for  any $f\in R[y,y',...,y^{( r )}]\,\widehat{}$ and $i\geq 0$:

$$
\d^{i+1} f(u)  =  \d \d^{i} f(u) 
$$
$$
= \frac{1}{p}((\d^{i}f)^{(\phi)}(u^p+p\d u,(\d u)^p+p\d^2 u,\dots )-(\d^i f)(u,\d u,\dots )^p)
$$
$$
\equiv   \frac{\partial \d^{i} f}{\partial p}(\nabla(u)^p)+ \sum_{j=0}^r \frac{\partial (\d^i f)^{(\phi)}}{\partial y^{(j)}} (\nabla(u)^p) \d^{j+1} u\ \roman{mod}\ p
$$
$$
 \equiv   \frac{\partial \d^{i} f}{\partial p}(\nabla(u)^p)+ \langle \, d [(\d^{i}f)^{(\phi)}] (\nabla(u)^p), 
\d \nabla(u)\cdot \left( \frac{\partial}{\partial \nabla(y)}\right)^t
\rangle\   \roman{mod}\ p. 
$$
It is known ([Bu05], Lemma 3.20), that if $b\in R$ then $b=0$ if and only if 
$$
\d^i b\equiv 0\ \roman{mod} \ p\ \ \roman{for\ all}\ i\geq 0.
$$

Combining the above facts we get:
\medskip

{\bf 5.2. Proposition.} {\it 
An element $u\in R$ is a solution to $f(u,\d u,...,\d^r u)=0$ if and only if the following hold:
\smallskip
a)  $f(\nabla(u))  \equiv   0 \ \roman{mod} \ p$,
\smallskip
b) $\dfrac{\partial (\d^{i} f)}{\partial p}(\nabla(u)^p)  \equiv   - \langle\, d [(\d^{i}f)^{(\phi)}] (\nabla(u)^p), 
\d \nabla(u)\cdot \left( \dfrac{\partial}{\partial \nabla(y)}\right)^t
\rangle \ \roman{mod} \ p,\ \ i\geq 0.$
\smallskip
 Moreover, $u$ is a solution to the equation $\d f(u,\d u,...,\d^r u)=0$ if and only if b) above holds.}

\smallskip
The above discussion can be obviously generalized to the case when $y$ is a tuple of variables
and $f$ is a tuple of equations. It shows that the equation $f(u,\d u,...,\d^r u)=0$ is controlled by a system of ``Pfaffian" congruences involving the $1$-forms
$$d((\d^i f)^{(\phi)})=d(\d^i (f^{(\phi)})),\ \ i\geq 0.$$
If $f$ has $\bold{Z}_p$-coefficients then the above forms are, of course,
$$d(\d^i f),\ \ i\geq 0.$$
In what follows we analyse such forms relevant for $PVI$ equations.

\medskip

{\bf 5.3. The forms $\delta^i(\psi-\rho)$.} We start by defining inductively certain universal $\d$-polynomials.
\smallskip
Consider a family of commuting free variables $z= (z_0,z_1, \dots ,z_r)$  and let $w$ be another variable. 
As in the jet theory in [Bu05], denote by
$z',z'', \dots $, resp.  $w',w'',...$  new (families of) independent  variables indexed
additionally by the formal order of derivative. Consider the ring
$\bold{Z}[z,z',z'' ,\dots ,w,w',w'' , \dots],$
equipped with the tautological $p$--derivation $\delta z_i^{(k)}= z_i^{(k+1)}$ etc., and hence with the lift of Frobenius $\phi(F)=F^p+p\d F$.  Define the elements $A_{m,i}=A_{m,i}(z,z',z'',\dots , w,w',w'\dots )$ of this ring ($m\geq 0$, $i=0,...,m+r$)  by induction:
$$
A_{0,i} :=  z^{(i)}, \quad  i=0,...,r
$$
$$
A_{m+1,0} :=  -(w^{(m)})^{p-1} A_{m,0}, \quad m\geq 0 \eqno(5.1)
$$
$$
A_{m+1,i} :=  \phi(A_{m,i-1})-(w^{(m)})^{p-1}A_{m,i}, \  i=1,...,m+r,\ \ m\geq 0
$$
$$
A_{m+1,m+r+1}  : =  \phi(A_{m,m+r}), \  m\geq 0.
$$

\medskip
It is easy then to check that
$$
A_{m,0} =  (-1)^{m-1}(ww'...w^{(m-1)})^{p-1}z^{(0)},\ \ m\geq 1,
$$
$$
A_{m,m+r} =  \phi^m(z^{(r)}),\ \ \ m\geq 0,
$$
$$
A_{m,m+r-1}  \equiv  \phi^m(z^{(r-1)})\ \roman{mod}\  (z^{(r)},\phi(z^{(r)}), \dots ,\phi^{m-1}(z^{(r)})),\ \ m\geq 0
\eqno(5.2)
$$
$$
A_{m,m+i}  \equiv  \phi^m(z^{(i)})\ \roman{mod}\  (w,w', \dots ,w^{(m-1)}),\ \ m\geq 1,\ \ i=0, \dots ,r,
$$
$$
A_{m,i}  \equiv  0\ \roman{mod}\  (w,w', \dots ,w^{(m-1)}),\ \  i=0, \dots ,m-1.
$$
\smallskip

In the following statement $Y$ is an affine smooth curve over $R$  and  $f\in \Cal{O}(J^r(Y))$. Then by [Bu05], 
there exist $a_0,...,a_r\in \Cal{O}(J^r(Y))$ such that
$$
df=\sum_{i=0}^r a_i\omega^{(i)}.
$$
Put also $a=(a_1, \dots , a_r)$. Then we have:
\medskip
{\bf 5.3.1. Proposition.} {\it Put $a_{m,i}=A_{m,i}(a,\d a, \d^2 a, \dots ,f,\d f, \d^2 f,...)$. Then
$$
d (\d^m f)=\sum_{i=0}^{m+r} a_{m,i} \omega^{(i)}.
$$
}
{\it Proof}. We proceed by induction on $m\geq 0$. The case $m=0$ is trivial.  The passage from $m$ to $m+1$
runs as follows:
$$
d(\d^{m+1} f)  =  d(\d \d^m f)
 =  d\left(\frac{\phi(\d^mf)-(\d^mf)^p}{p}\right)
=  \frac{\phi^*}{p}(d(\d^m f))-(\d^m f)^{p-1} d(\d^m f)
$$
$$
 =  \frac{\phi^*}{p}(\sum_{i=0}^{m+r}a_{m,i}\omega^{(i)})-(\d^m f)^{p-1} (\sum_{i=0}^{m+r} a_{m,i} \omega^{(i)})
$$
$$
=  (\sum_{i=0}^{m+r}\phi(a_{m,i})\omega^{(i+1)})-(\d^m f)^{p-1} (\sum_{i=0}^{m+r} a_{m,i} \omega^{(i)})
=  \sum_{i=0}^{m+1+r} a_{m+1,i}\omega^{(i)}.
$$
This ends the proof.

\medskip
Now we apply this to the case when  $Y$ is an affine open subset of the elliptic curve $E$ over  $\bold{Z}_p$
without lift of Frobenius. Denote by $a_p\in \bold{Z}$ the trace of Frobenius on the reduction
$E\,\roman{mod}\,p.$

\medskip
{\bf 5.3.2. Corollary.} {\it Let $f=\psi-\rho$, $\rho=\phi(r)$, where $r\in \Cal{O}(Y)$,
$\psi$ the canonical $\d$-character of order $2$  attached to $\omega$. Let $d(\d^m f)=\sum_{i=0}^{m+2} a_{m,i} \omega^{(i)}$, $m\geq 0$, $a_{m,i}\in \Cal{O}(J^{m+2}(Y))$. So 
$$
a_{0,2}=p,\ a_{0,1}=-(a_p+p\phi(\frac{dr}{\omega})),\ a_{0,0}= 1.
$$
 Then, for $m\geq 1$:
\smallskip
1) $a_{m,m+2}=p$;

2) $a_{m,m+1}\equiv -a_p\, \roman{mod}\, p$;

3) $a_{m,m+1}\equiv -(a_p+p\phi^{m+1}(\frac{dr}{\omega}))\, \roman{mod}\, (f,\d f, \dots ,\d^{m-1}f)$;

4) $a_{m,m}\equiv 1\, \roman{mod}\, (f,\d f, \dots ,\d^{m-1}f)$;

5) $a_{m,i}\equiv 0\,\roman{mod}\,(f,\d f, \dots ,\d^{m-1}f)$ for $i=0, \dots , m-1$;

6) $a_{m,0}  =  (-1)^{m-1}(f \cdot \d f \cdot \cdot \cdot \d^{m-1} f)^{p-1}$, for $m\geq 1$.

\smallskip
In particular if $E$ has ordinary reduction then $a_{m,m+1}$ is invertible in $\Cal{O}(J^{m+2}(Y))$.}

\medskip

One can prove a similar statement for $\psi-r$ in place of $\psi-\phi(r)$.

\medskip

{\bf 5.3.3. Corollary.} {\it
In the same situation, let $f=\psi-\rho$, where $\rho=r$ or $\rho=\phi ( r )$.   Let $\psi$ be the canonical $\d$-character of order $2$  attached to $\omega$.
\smallskip
Denote by $Z^{m+2}=Z^{m+2}(f,\d f,...,\d^m f)$  the closed formal subscheme of $J^{m+2}(Y)$ defined by the ideal generated by $f$, $\d f$, ..., $\d^m f$. Let $\pi_m:Z^{m+2}\to J^1(Y)$ be the canonical projection.
\smallskip
Then the map
$$
\pi_m^*\Omega_{J^1(Y)}\to \Omega'_{Z^{m+2}}
$$
is injective with  cokernel annihilated by $p^{m+1}$. In particular, for any $m\geq 0$,  the system $\{f,\d f,...,\d^m f\}\subset \Cal{O}(J^{m+2}(Y))$ defines a generalised $\d$--flow on $J^1(Y)$. If moreover $E$ has ordinary reduction then the above cokernel is a cyclic $\Cal{O}(Z^{m+2})$-module generated
by the class of $\omega^{(m+2)}$.}

\medskip

{\bf Proof}. Consider the case $\rho=\phi(r)$; a similar argument holds for $\rho=r$.
Recall that $\Omega_{J^n(Y)}$ is a free $\Cal{O}(J^n(Y))$-module generated by 
$\omega^{(0)},...,\omega^{(n)}$. Also, by definition,
$$\cO(Z^{m+2})=\frac{\cO(J^{m+2}(Y))}{(f,\d f,...,\d^m f)}$$
and
$$\Omega'_{Z^{m+2}}=\frac{\cO(J^{m+2}(Y))\omega^{(0)}\oplus...\oplus \cO(J^{m+2}(Y))\omega^{(m+2)}}{\langle (\d^i f)\omega^{(j)}, d (\d^i f)\rangle},$$
where $\langle\ \ \rangle$ means $\cO(J^{m+2}(Y))$-linear span and $j=0,...,m+2$, $i=0,...,m$. By Corollary 5.3.2 we have:
$$
\Omega'_{Z^{m+2}}=\frac{\cO(Z^{m+2})\omega^{(0)}\oplus ... \oplus \cO(Z^{m+2})\omega^{(m+2)}}{\langle p\omega^{(2+i)}-(a_p +p\phi^{i+1}(\frac{dr}{\omega}))\omega^{(1+i)}+\omega^{(i)}\rangle} 
$$
where $i=0,...,m$ and $\langle\ \ \rangle$ means here $\cO(Z^{m+2})$-linear span. Note that $\pi^*_m\Omega_{J^1(Y)}$ is a free $\cO(J^{m+2}(Y))$-module 
with basis $\omega^{(0)}$, $\omega^{(1)}$. So in  order to prove that the map
$$
\pi^*_m\Omega_{J^1(Y)}\to \Omega'_{Z^{m+2}}
\eqno(5.3)
$$
is injective we need to check that no $\cO(Z^{m+2})$-linear combination of 
$\omega^{(0)}$, $\omega^{(1)}$ can be a $\cO(Z^{m+2})$-linear combination
of elements $p\omega^{(2+i)}-(a_p+p\phi^{i+1}(\frac{dr}{\omega})) \omega^{(1+i)}+\omega^{(i)}$ which is clear. For $j=2,...,m+2$ let 
$$\overline{\omega^{(j)}}\in \frac{\Omega'_{Z^{m+2}}}{\langle \omega^{(0)},\omega^{(1)}\rangle}$$
be the image of $\omega^{(j)}$. Clearly $p\overline{\omega^{(2)}}=0$ hence 
$p^2\overline{\omega^{(3)}}=0$, etc. We conclude that $p^{m+1}$ annihilates the cokernel of (5.3).

\smallskip
 Finally assume $E$ has ordinary reduction. We want to show that the above cokernel  is generated by $\overline{\omega^{(m+2)}}$. It is enough to show that for all $j=1,...,m+2$ we have  
$$
\overline{\omega^{(j-1)}}\in p\,\cO(Z^{m+2})\,\overline{\omega^{(j)}}.
$$
We proceed by induction on $j$. For $j=1$ this is clear. Now assume the above is true  for some $1\leq j<m+2$. We have 
$$
p\overline{\omega^{(j+1)}}-(a_p+p\phi^j(\frac{dr}{\omega})) \overline{\omega^{(j)}}+\overline{\omega^{(j-1)}}=0.
$$
By induction $\overline{\omega^{(j-1)}}=p\,c\,\overline{\omega^{(j)}}$ for some $c\in \cO(Z^{m+2})$. Since $a_p$ is invertible in $R$ it follows that $a_p-pc(1+\phi^{j-1}(dr/\omega))$ is invertible hence
$$\overline{\omega^{(j)}}=p(a_p+p\phi^j(\frac{dr}{\omega})-pc)^{-1}\overline{\omega^{(j+1)}},$$
which ends the proof.

\medskip

{\bf 5.4. Remark.}
Let $Y$ be a smooth affine scheme over $R$ and let $\Cal{F}\subset \cO(J^r(Y))$ be a system of arithmetic differential equations. Consider the set 
$$
\{\Cal{F}, \d \Cal{F}, \dots ,\Cal{F} \} \subset \cO(J^{r+m}(Y))
$$
which can be referred to as the $m$--th prolongation of $\Cal{F}$. Let
 $$
 Z^{r+m}:=Z^{r+m}(\Cal{F}, \d \Cal{F} , \dots ,\d^m \Cal{F})\subset J^{r+m}(Y)
 $$ 
 be the closed subscheme defined by this prolongation and let 
$$
Z^{\infty}=Z^{\infty}(\d^{\infty}\Cal{F})
$$
be defined by
$$
Z^{\infty}=\lim_{\leftarrow} Z^{r+m}
$$ 
be the projective limit (defined as the $Spf$ of the $p$--adic completion of the inductive limit of $\cO(Z^{r+m})$ as $m$ varies; this is a generally non-Noetherian formal scheme). Then $Z^{\infty}$ is a closed horizontal formal subscheme of  
$$
J^{\infty}(Y):=\lim_{\leftarrow} J^n(Y);
$$
by horizontal we mean here that its ideal is sent into itself by $\d$. So there is an induced $p$--derivation 
$\d_{Z^{\infty}}$ on $\cO(Z^{\infty})$. 
\smallskip
If $Z^{\infty}$ happens to be a smooth formal scheme then $\d_{Z^{\infty}}$ defines a (genuine) $\d$-flow
$\Cal{E}(\d_{Z^{\infty}})\subset \cO(J^1(Z^{\infty}))$ on $Z^{\infty}$ (in the sense of our previous definition). In case $Z^{\infty}$ is not necessarily a smooth formal scheme one can still define the affine formal scheme 
$J^1(Z^{\infty})$ and the $\d$--flow
$\Cal{E}(\d_{Z^{\infty}})$ on $Z^{\infty}$ by copying the definitions from the smooth case.
\smallskip
The following problem needs to be investigated. Assume that $\Cal{F}\subset \cO(J^2(Y))$ defines a $\d$-flow 
$\Cal{E}(\d_{J^1(Y)})\subset \cO(J^1(J^1(Y)))$ on $J^1(Y)$. 
Consider the formal scheme $Z^{\infty}=Z^{\infty}(\d^{\infty} \Cal{F})$ as above.
Prove that $Z^{\infty}$ is naturally isomorphic to $J^1(Y)$ and $\Cal{E}(\d_{J^1(Y)})$ is naturally identified with $\Cal{E}(\d_{Z^{\infty}})$. This would show the naturality of the above definitions.

\vskip1cm

\centerline{\bf References}

\medskip

[BlSch12] M.~Blickle, K.~Schwede. {\it $p^{-1}$--linear maps in algebra and geometry.}
arXiv:1205.4577
\smallskip
[Bo09]  J.~Borger. {\it Lambda--rings and the field with one element.}  arXiv:0906.3146 
    
\smallskip
[Bu95] A.~Buium. {\it Differential characters of abelian varieties over $p$--adic fields.} Inv.~Math.,
vol. 122 (1995), 309--340.

\smallskip

[Bu97] A.~Buium. {\it Differential characters and  characteristic polynomial of Frobenius.} 
J. reine u. angew. Math. 485 (1997), 209--219.

\smallskip
[Bu05] A.~Buium. {\it Arithmetic Differential Equations.} Math. Surveys and Monographs,
118, MS, Providence RI, 2005. xxxii+310 pp.

\smallskip 
[Ma96] Yu.~Manin. {\it Sixth Painlev\'e equation, universal elliptic curve,
and mirror of $\bold{P}^2$.}  In: geometry of Differential
Equations, ed. by A.~Khovanskii, A.~Varchenko, V.~Vassiliev.
Amer. Math. Soc.
Transl. (2), vol. 186  (1998), 131--151. Preprint alg--geom/9605010.
\smallskip
[Ma08] Yu.~Manin. {\it Cyclotomy and analytic geometry over $F_1$.}  In: Quanta of Maths. Conference in honour of Alain Connes.
Clay Math. Proceedings, vol. 11 (2010), 385--408.
 Preprint math.AG/0809.2716.
\smallskip

[Ta01] K.~Takasaki. {\it Painlev\'e--Calogero correspondence revisited.} Journ.~Math.~Phys., vol.~42, Nr 3 
(2001), 1443--1473.

\enddocument